\documentclass[a4paper,12pt]{article}

\usepackage{amsmath,amssymb,amsfonts}
\usepackage{amsthm,epsfig,epstopdf,url,array}
\usepackage{mathtools}
\usepackage{enumerate}
\usepackage{graphicx}
\theoremstyle{plain}
\newtheorem{thm}{Theorem}[section]
\newtheorem{lem}[thm]{Lemma}
\newtheorem{prop}[thm]{Proposition}
\newtheorem*{cor}{Corollary}

\theoremstyle{definition}
\newtheorem{defn}{Definition}[section]

\theoremstyle{remark}
\newtheorem*{rem}{Remark}

\usepackage[a4paper,includeheadfoot,margin=2.54cm]{geometry}

\begin{document}

\title{Categorical study for Algebras of lattice-valued logic and lattice-valued modal logic}
\author{Kumar Sankar Ray\thanks{E-mail:ksr@isical.ac.in} \and  Litan Kumar Das\thanks{E-mail:ld06iitkgp@gmail.com}}

 \maketitle

{\bf Abstract.} 
                    The paper explores categorical interconnections between lattice-valued Relational systems and algebras of Fitting's lattice-valued modal logic. We define lattice-valued boolean systems, and then we study co-adjointness, adjointness of functors. As a result, we get a duality for algebras of lattice-valued logic. Following this duality results, we establish a duality for algebras of lattice-valued modal logic.\\

\section*{Introduction}
Vicker's in \cite{vickers1996topology} introduced the concept of topological systems in the work of topology via logic, which was further considered in \cite{vickers1999topology}. Topological systems is a mathematical structure like as $(X,A,\models)$, where $X$ is a non-empty set, $A$ is a frame, a complete distributive lattice and $\models$ is a satisfaction relation on $X\times A$. This relation $\models$ satisfies both join and finite meet interchange laws.\\
The authors in \cite{denniston2009lattice,denniston2012interweaving, denniston2009functorial} demonstrated the idea of lattice-valued topological systems and considered the category $\ell$-$\textbf{TopSys}$ from lattice-valued topological systems. They also explored categorical relation with the systems and spaces. Besides, in \cite{solovyov2010variable} variable basis-topological systems was taken place as a further generalization of lattice-valued topological systems, which was again considered in another point of view(see \cite{solovyov2011localification}).\\
In \cite{fitting1991many} Fitting revealed the idea of $\ell$-valued logic and $\ell$-valued modal logic for a finite distributive lattice $\ell$, where $\ell$ is endowed with the truth constants. Several studies have been done in various aspects of Fitting's logics (see \cite{fitting1991many,fitting1992many,fitting1995tableaus,eleftheriou2005frame,koutras2003catalog}). Maruyama \cite{maruyama2009algebraic} defined the class of $\ell$-\textbf{VL}-algebras and the class of $\ell$-\textbf{ML}-algebras as an algebraic structure of Fitting's $\ell$-valued logic and $\ell$-valued modal logic respectively. Consequently, a duality developed for $\ell$-\textbf{VL}-algebras and the category $\ell$-\textbf{BS} in \cite{maruyama2011dualities}, which can be conceived as strong duality according to the theory of natural dualities \cite{clark1998natural}. Following the duality for $\ell$-\textbf{VL}-algebras, he also cultivated J$\acute{o}$sson-Tarski duality(see\cite{chagrov1997modal, hansoul1983duality}) for $\ell$-\textbf{ML}-alegebras.\\
While studying \cite{denniston2012interweaving}, we raised a question whether there exists systems which are categorically connected with $\ell$-\textbf{VL}-algebras and $\ell$-\textbf{ML}-algebras. 
Our objective here is to define such systems and show this systems are categorically equivalent with the spaces(c.f.\cite{maruyama2011dualities}). 
As a result, we shall establish a duality for $\ell$-\textbf{VL}-algebras and $\ell$-\textbf{ML}-algebras. \\
The paper is arranged as follows. \\
In Section [\ref{sec:1}], we recall some basic notions associated with this work. We define a category $\ell$-\textbf{BSYM} from the concept of $\ell$-boolean systems in Section [\ref{sec:2}] and establish categorical relationships with $\ell$-Boolean space $\ell$-\textbf{BS} and $\ell$-\textbf{VL}-algebras. In Section [\ref{sec:3}], we introduce the concept of $\ell$-relational systems and define a category  $\ell$-\textbf{RSYM}. Henceforth we shall show duality between $\ell$-\textbf{ML}-algebras and $\ell$-relational space $\ell$-\textbf{RS}. The paper is concluded in Section [\ref{sec:4}]

\section{Preliminaries}
\label{sec:1}
For category theory we refer to \cite{adamek1990h, mac2013categories}. To read our paper easily, we mention here some crucial concepts.\\
\begin{defn}
Let $f: \textbf{G}\rightarrow \textbf{H}$ be a functor, and take $H$ be a $\textbf{H}$-object.
\begin{enumerate}
\item A $f$-structured arrow with domain $H$ is a pair $(g,G)$ consisting of an $\textbf{G}$-object $G$ and a \textbf{H}-morphism $g: H\rightarrow f(G)$.
\item A $f$-structured arrow with domain $H$ is called $f$-universal arrow for $H$ provided that for each $f$-structured arrow $(g', G')$ with domain $H$ there exists a unique \textbf{G}-morphism 
$\tilde{g}: G\rightarrow G'$ with $g'=f(\tilde{g})\circ g$, in otherwords the triangle commutes.
\begin{figure}
\begin{center}
\includegraphics[width=90mm]{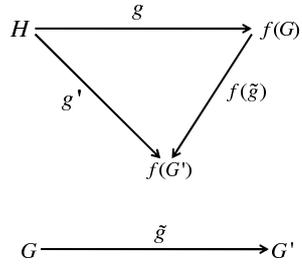}
\caption{Diagram of the Universal arrow}
\label{fig:univ}
\end{center}
\end{figure}
\item A $f$-costructured arrow with codomain $H$ is a pair $(G, g)$ consisting of a \textbf{G}-object $G$ and a \textbf{H}-morphism $g: f(G)\rightarrow H$
\item A $f$-costructured arrow $(G, g)$ with codomain $H$ is called $f$-couniversal for $H$ provided that for each $f$-costructured arrow $(G',g')$ with codomain $H$ there exists a unique \textbf{G}-morphism $\tilde{f}: G'\rightarrow G$ with $g'=g\circ f(\tilde{g})$.
\end{enumerate}
\end{defn}

\begin{defn}[adjoint]
A functor $f$: \textbf{G}$\rightarrow$ \textbf{H} is said to be adjoint if for every \textbf{H}-object $H$ there exists a $f$-universal arrow with domain $H$. Consequently, there exists unique natural transformation (called the unit) $\eta_H: id_H\rightarrow ff_1$, where $id_H$ is an identity morphism from $H$ to $H$ and $f_1: H\rightarrow G$ is a functor, is a $f$-universal arrow. Diagram of unit is shown in Fig.\ref{fig:unitsimple})
\begin{figure}
\begin{center}
\includegraphics[width=100mm]{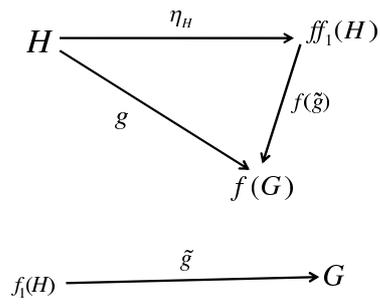}
\caption{Diagram of Unit}
\label{fig:unitsimple}
\end{center}
\end{figure}
\end{defn}
\begin{defn}[coadjoint]
A functor $f$: \textbf{G}$\rightarrow$ \textbf{H} is said to be co-adjoint if for every \textbf{H}-object $H$ there exists a $f$-couniversal arrow with codomain $H$. As a result, there exists a natural transformation $\xi_G: f_1\circ f(G)\rightarrow id_G(G)$, where $id_G$ is an identity morphism from $G$ to $G$, and $f_1: H\rightarrow G$ is a functor such that for a given morphism $g: f(H)\rightarrow G$, there is a unique \textbf{H}-morphism $\tilde{g}: H\rightarrow f(G)$ such that the triangle (see Fig.\ref{fig:counitsimp}) commutes, in otherwords $g=f_1(\tilde{g})\circ \xi_G$ . Here $\xi$ is called the co-unit of the adjunction. Diagram of counit is shown in Fig.\ref{fig:counitsimp}
\begin{figure}
\begin{center}
\includegraphics[width=100mm]{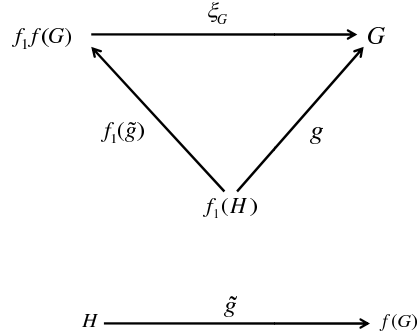}
\caption{Diagram of Co-unit}
\label{fig:counitsimp}
\end{center}
\end{figure}

\end{defn}

\begin{defn}(\cite{vickers1996topology})[Frame]
A poset(partially ordered set) is said to be frame iff
\begin{enumerate}
\item every subset has a join
\item every finite subset has a meet
\item binary meets distribute over joins:\\
$x\wedge\bigvee=\bigvee \{x\wedge y: y\in Y\}$
\end{enumerate}
\end{defn}

\section{$\ell$-\textbf{VL}-algebras, $\ell$-Boolean Systems, $\ell$-Boolean Space and their Categorical interconnections} 
\label{sec:2}
\subsection{$\ell$-\textbf{VL}-algebras}
$\ell$-valued logic $\ell$-\textbf{VL} is a many-valued logic and $\ell$-\textbf{VL}-algebras defined in \cite{maruyama2009algebraic} as an algebraic semantics for $\ell$-\textbf{VL}, which are both sound and complete.
\begin{defn}(\cite{maruyama2009algebraic})
$\ell$-\textbf{VL}-algebra homomorphism is a function between two $\ell$-\textbf{VL}-algebras which preserves the operations $(\vee, \wedge, \rightarrow, T_r(r \in \ell), 0,1)$.

\end{defn}
\begin{defn}(\cite{maruyama2011dualities})
Let A be an $\ell$-\textbf{VL}-algebra. A non-empty subset $F$ of $A$ is called an $\ell$-filter iff $F$ is a filter of lattices which is closed under $T_1$. Let $P$ be a proper $\ell$-filter of $A$.\\
\begin{enumerate}
\item P is a prime $\ell$-filter of $A$ iff for any $r\in \ell$, $T_r(x\vee y)\in P$, then there exist $r_1,r_2\in \ell$ with $r_1\vee r_2=r$ such that $T_{r_1}(x) \in P$ and $T_{r_2}(y)\in P$.
\item P is an ultra $\ell$-filter of $A$ iff $\forall r\in A$ $\exists r\in \ell$ ,$T_r(x)\in P$.
\item $P$ is a maximal $\ell$-filter iff $P$ is maximal with respect to inclusion.
\end{enumerate}
\end{defn}
\begin{prop}(\cite{maruyama2011dualities}
\label{propo_1}
\begin{enumerate}
\item Let $A$ be an $\ell$-$\textbf{VL}$-algebra. For any two distinct members $x,y$ of $A$ ,there exist $r\in \ell$ and a prime $\ell$-filter $P$ of $A$ such that $T_r(x)\in P$ and $T_r(y)\notin P$.
\item For a prime $\ell$-filter $P$ of an $\ell$-$\textbf{VL}$-algebra $A$, define $v_P:A \rightarrow \ell$ by $v_P(x)=q \Leftrightarrow T_q(x)\in P$. Then, $v_P$ is a homomorphism of $\ell$-\textbf{VL}-algebras.
\item Let $A$ be an $\ell$-\textbf{VL}-algebra. A bijective mapping exists from the set of all prime $\ell$-filters of $A$ to the set of all homomorphisms from $A$ to $\ell$.

\end{enumerate}
\end{prop}

The spectrum of an $\ell$-\textbf{VL}-algebra $A$ is designated by $Spec_{\ell}(A)$, and is defined as follows.
\begin{defn}(\cite{maruyama2011dualities})\label{spec}
Let $A$ be an $\ell$-\textbf{VL} algebra. For a subalgebra $\ell_1$ of $\ell$, $Spec_{\ell_1}(A)=$$\{\ell$-\textbf{VL}-algebras homomorphism $f:A \rightarrow \ell_1\}$.
\end{defn}
\subsection{Categories and their Functorial relationships}
\textbf{$\ell$-\textbf{Boolean systems}}
\begin{defn}
\label{bsym}
An $\ell$-Boolean system is a triple $(X,A,\models)$, where $X$ is a non-empty set, $A$ is an $\ell$-\textbf{VL}-algebra and $\models$ is an $\ell$-valued satisfaction relation on $(X,A)$, i.e., $\models : X\times A \rightarrow \ell$ is a mapping such that 
\begin{enumerate}
\item if $A_1$ is a subset of $A$, then $\displaystyle \models(x,\bigvee_{a \in A_1})=\bigvee_{a\in A_1}\models(x,a)$, \\
\hspace{0.5in} if $A_1$ is a finite subset of $A$, then $\displaystyle\models(x,\bigwedge_{a\in A_1})=\bigwedge_{a\in A_1}\models(x,a)$
\item if $x_1\neq x_2$ in $x$ then $\displaystyle\models(x,a_1) \neq \models(x_2,a)$, for some $a\in A$.
\item $\models(x,a\rightarrow b)=\models(x,a)\rightarrow \models(x,b)$.
\item $\models(x,T_r(a))=T_r(\models(x,a))$, $a\in A$ and $r \in \ell$.
\item $\models(x,0)=0$, $\models(x,1)=1$. Here 1, 0 are respectively the top-element and bottom-element.

\end{enumerate}

\end{defn}

\begin{defn}
\label{cbsym}
We define a category $\ell$-\textbf{BSYM} as follows. \\
(1) Objects are all $\ell$-Boolean systems $(X,A,\models)$.\\
(2) Arrows are all continuous maps $(\psi_1,\psi_2):(X,A,\models_1)\rightarrow (Y,B,\models_2)$, where 
\begin{center}
\begin{itemize}
\item $\psi_1:X\rightarrow Y$ is a set map.
\item $\psi_2:B\rightarrow A$ is a homomorphism of $\ell$-\textbf{VL}-algebras.
\item $\models_1(x,\psi_2(b))=\models_2(\psi_1(x),b)$
\end{itemize}
\end{center}
(3) For each object $P=(X, A, \models)$ ,the identity arrow $I_P: P\rightarrow P$ is the pair $(I', I'')$ such that 
\begin{center}
$I':X\rightarrow X$ \\
$I'':A\rightarrow A$.
\end{center}
(4) For a given $\ell$-Boolean systems $P=(X, A, \models_1)$, $Q=(Y, B, \models_2)$, $R=(Z, C, \models_3)$ let $(\psi_1, \psi_2): P\rightarrow Q$ and $(\phi_1, \phi_2): Q\rightarrow R$ be continuous maps. Composition defined as $ (\phi_1, \phi_2)\circ(\psi_1, \psi_2): P\rightarrow R$ such that 
\begin{center}
$\phi_1\circ\psi_1:X\rightarrow Z$ \\
$\psi_2\circ\phi_2:C\rightarrow A$.
\end{center}

\end{defn}
\begin{defn}We introduce the notion of \textit{extent} in an $\ell$-Boolean system $P=(X,A,\models)$.
For each member $a$ in $A$, its \textit{extent} in $P$ is a function $ext_{\ell}(a): X\rightarrow \ell$ such that $ext_{\ell}(a)(x)=\models(x, a)$.\\
$ext_{\ell}(A)=\{ext_{\ell}(a): a\in A\}$.On the set $ext_{\ell}(A)$, the operations $(\wedge, \vee, \rightarrow, T_r, 0, 1)$ are defined pointwise. 
\end{defn}
A zero-dimensional compact Hausdorff space is said to be Boolean space(see\cite{johnstone1986stone}).
\begin{cor}
$(X,ext_{\ell}(A))$ is a Boolean space.
\end{cor}

Designated by $\ell$-VA, the category of $\ell$-\textbf{VL}-algebras and homomorphism of $\ell$-\textbf{VL}-algebras.\\
The author in (\cite{maruyama2011dualities}) defined a category $\ell$-BS. Here we note that $\ell$ and its subalgebras denoted by $\it{Subalg}(\ell)$ are taken with the discrete topology.Then $(\ell, \phi_{\ell})$ is an object in $\ell$-BS.
Where $\phi_{\ell}:Subalg(\ell)\rightarrow Subsp(\ell)$ such that $\phi_{\ell}(S)=S$,$S\in Subalg(\ell)$.

\begin{defn}(\cite{maruyama2011dualities})
Let $(\mathcal{O}, \phi)$ be an object in the category $\ell$-BS. Then $Cont(\mathcal{O}, \phi)$ is defined as the set of all subspace-preserving continuous functions from $(\mathcal{O}, \phi)$ to $(\ell, \beta)$. $Cont(\mathcal{O}, \phi)$ is endowed with the operations $(\wedge, \vee,\rightarrow, T_r(r\in \ell), 0, 1)$ defined pointwise.

\end{defn}

We shall now study functorial relationship among the categories $\ell$-\textbf{BSYM}, $\ell$-BS and $\ell$-VA.
\begin{defn}
we define a functor $Ext_{\ell}$: $\ell$-\textbf{BSYM}$\rightarrow \ell$-BS as follows.

For an object $(X,A,\models)$ in $\ell$-\textbf{BSYM}, define $Ext_{\ell}(X,A,\models)=(X,ext_{\ell}(A),\phi)$. Where $\phi:Subalg(\ell)\rightarrow Subsp(S)$ such that $\phi(K)=(X,ext_K(A))$, $K\in Subalg(\ell)$.\\
For an arrow $(f,\phi):(X,A,\models_1)\rightarrow (Y,B,\models_2)$, define $Ext_{\ell}(f,\phi)=f$.\\ Where $f:(X,ext_{\ell}(A),\alpha_1)\rightarrow (Y,ext_{\ell}(B),\alpha_2)$ is a subspace preserving continuous map.

\end{defn}
\begin{rem}
Here $f=Ext_{\ell}(f,\phi):(X,ext_{\ell}(A))\rightarrow (Y,ext_{\ell}(B))$ is continuous, since $f^{-1}(ext_{\ell}(b)(x)=ext_{\ell}(b)f(x)=\models_2(f(x),b)=\models_1(x,\phi(b))=ext_{\ell}(\phi(b))(x)$. Therefore $f^{-1}(ext_{\ell}(b))=ext_{\ell}(\phi(b)) \in ext_{\ell}(A)$.
\end{rem} 
\begin{lem}
Let $(R,\alpha)$ be an object in $\ell$-BS. Then $(R,Cont(R,\alpha),\models)$ is an object in $\ell$-\textbf{BSYM}.
\begin{proof}
Here $Cont(R,\alpha)$ is an $\ell$-{\textbf{VL}}-algebra and $\models(r,v)=v(r)$. Now we verify that $Cont(R,\alpha)$ is an $\ell$-Boolean system.
\begin{enumerate}[(i)]
\item For a subet $A_1$ of $A$ , $\models(r,\vee_{a\in A_1}u_a)=(\vee_{a\in A_1}u_a)(r)=\vee_{a\in A_1}u_a(r)=\vee_{a\in A_1}\models(r,u_a)$.\\
For any $u_1, u_2\in R$, $\models(r,u_1\wedge u_2)=(u_1\wedge u_2)(r)=u_1(r)\wedge u_2(r)=\models(r,u_1)\wedge \models(r,u_2)$.
\item As R is zero-dimenssional and Hausdorff space and hence Kolmogorov, for $r_1\neq r_2$ in R there exists $v\in Cont(R,\alpha)$ such that $v(r_1)\neq v(r_2)$. So $\models(r_1,v)\neq \models(r_2,v)$, for some $v\in Cont(R,\alpha)$.
\item $T_{\ell}(\models(r,v))=T_{\ell}(v(r))=T_{\ell}(v)(r)=\models(r,T_{\ell}(v))$.
\item $\models(r,v\rightarrow u)=(v\rightarrow u)(r)=v(r)\rightarrow u(r)$.                                                                                                                                                                                                                                                                                                                                                                                                                                                                                                                                                                                                                                                                                                                                                                                                                                                                                                                                                                                                                                                                                                                                                                                                                                                                                                                                                                                                                                                                                                                                                                                                                                                                                                                                                                                                                                                                                                                                                                                                                                                                                                                                                                                                                                                                                                                                                                                                                                                                                                                                                                                                                                                                                                                                                                                                                                                                                                                                                                                                                                                                                                                                                                                                                                                                                                                                                                                                                                                                                                                                                                                                                                                                      

\end{enumerate}
\end{proof}
\end{lem}
\begin{lem}
If $f$ is an arrow in $\ell$-BS, then $(f, f^{-1})$ is an $\ell$-continuous map in $\ell$-\textbf{BSYM}.
\begin{proof}
Here $f:S\rightarrow T$ is a set function and $f^{-1}:Cont(T,\alpha_2)\rightarrow Cont(S,\alpha_1)$ is an $\ell$-VL-algebra homomorphism. $f^{-1}(v)=v\circ f$ for $v\in Cont(T,\alpha_2)$. Now $\models_2(f(s),v)=v(f(s)=f^{-1}(v)(s)=\models_1(s,f^{-1}(v))$. So $(f,f^{-1})$ is an $\ell$-continuous map in $\ell$-\textbf{BSYM}.
\end{proof}
\end{lem}

\begin{defn} We define a functor $P$: $\ell$-BS $\rightarrow \ell$-\textbf{BSYM}.
\begin{enumerate}[(i)]
\item P acts on an object in $\ell$-BS as $P(S,\alpha)=(S,Cont(S,\alpha),\models)$.
\item P acts on a morphism $f:(S,\alpha_1)\rightarrow (T,\alpha_2)$ in $\ell$-BS as $P(f)=(f,f^{-1}):(S,Cont(S,\alpha_1),\models_1)\rightarrow(T,Cont(T,\alpha_2),\models_2)$.\\
Where
\end{enumerate}
\begin{enumerate}[(i)]
\item $f:S\rightarrow T$, a set function.
\item A $\ell$-VL-algebra homomorphism $f^{-1}: Cont(T,\alpha_2)\rightarrow Cont(S,\alpha_1)$, defined as $f^{-1}(v)=v\circ f$.
\end{enumerate}
\end{defn}
\begin{defn}
Functor $Q$: $\ell$-$\textbf{BSYM}\rightarrow$ ($\ell$-VA)$^{op}$ acts on an object $(X,A,\models)$ in $\ell$-\textbf{BSYM} as $Q(X,A,\models)=A$, and on an arrow $(\psi_1,\psi_2): (X,A,\models_1)\rightarrow (Y,B,\models_2)$ in $\ell$-\textbf{BSYM} as $Q(\psi_1,\psi_2)=\psi_2: B\rightarrow A$ as $\psi_2$ is an $\ell$-VL-algebra homomorphism in ($\ell$-VA)$^{op}$.
\end{defn}
\begin{lem}
Let $A$ be an $\ell$-VL-algebra. Then $(Spec_{\ell}(A),A,\models_{(Spec_{\ell}(A)\times A)})$ is an object in $\ell$-\textbf{BSYM}.
\begin{proof}
Here $Spec_{\ell}(A)$ is a set, $A$ is a frame. Define $\models_{(Spec_{\ell}(A)\times A)}(v,a)=v(a)$. Now we verify the following \\
\begin{enumerate}[(i)]
\item $\displaystyle \models_{Spec_{\ell}(A)\times A}(v,\bigvee_{a\in \Gamma}a)=v(\bigvee_{a\in \Gamma}a)=\bigvee_{a\in \Gamma}v(a)=\bigvee_{a\in \Gamma}\models_{Spec_{\ell}(A)\times A}(v,a)$.                  
\item $\models_{(Spec_{\ell}(A)\times A)}(v,a_1\wedge a_2)=v(a_1\wedge a_2)=v(a_1)\wedge v(a_2)=\models_{(Spec_{\ell}(A)\times A)}(v,a_1)\wedge \models_{(Spec_{\ell}(A)\times A)}(v,a_2)$.
\item For $r\in \ell$, $T_r(\models_{(spec_{\ell}(A)\times A)}(v,a)=T_r(v(a))=v(T_r(a))$, and $\models_{(Spec_{\ell}(A)\times A)}(v,T_r(a))=v(T_r(a))$. 
Therefore, $\models_{(Spec_{\ell}(A)\times A)}(v,T_r(a))=T_r(\models_{(Spec_{\ell}(A)\times A)}(v,a))$. Others properties are checked easily.

\end{enumerate}
\end{proof}
\end{lem}
\begin{defn}
$R$ is a functor from ($\ell$-VA)$^{op}$ to $\ell$-\textbf{BSYM} defined as follows. \\
For an object $A$ in ($\ell$-VA)$^{op}$, $R(A)=(spec(A),A,\models_{(spec(A)\times A)})$ and on an arrow $f$ as $R(f)=(f^{-1}, f)$.
\end{defn}
By the above lemma it can be shown that $R$ is a functor.
\begin{thm}
$Ext_{\ell}$ is the co-adjoint to the functor $P$.
\begin{proof}
We shall prove the theorem establishing counit of the adjunction. Diagram of counit is shown in Fig. \ref{fig:coadjoint0}.
\begin{figure} 
\begin{center}
\includegraphics[width=100mm]{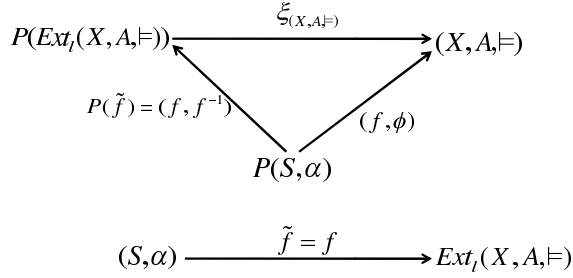}
\caption{Diagram of Counit}
\label{fig:coadjoint0}
\end{center}
\end{figure}
Here $P(S,\alpha)=(S,Cont(S,\alpha),\models)$ and $Ext_{\ell}(X,A,\models)=(X,ext_{\ell}(A),\alpha)$. Hence $P(Ext_{\ell}(X,A,\models))=P(X,ext_{\ell}(A),\alpha)
=(X,Cont(X,ext_{\ell}(A),\alpha),\models)$. Counit $\xi_{(X,A,\models)}=(id_X,ext_{\ell}): P(Ext_{\ell}(X,A,\models)\rightarrow (X,A,\models)$\\
 where 
\begin{enumerate}[(i)]
\item $id_X:X\rightarrow X$
\item $ext_{\ell}:A\rightarrow Cont(X,ext_{\ell}(A),\alpha)$
\end{enumerate}
$\mu$ is an open $\ell$ continuous map in $(\ell,\alpha_{\ell})$. $ext_{\ell}(a)\in Cont(X,ext_{\ell}(A),\alpha)$ since $ext_{\ell}(a)^{-1}(\mu)(x)=\mu\circ ext_{\ell}(a)(x)=\mu(\models(x,a))=\displaystyle\bigvee_{p\in \ell}(T_{\mu(p)}(ext_{\ell}(a)(x))=\bigvee_{p\in\ell} T_{\mu(p)}(\models(x,a))=\bigvee_{p\in \ell}(\models(x,T_{\mu(p)}(a))=\bigvee_{p\in\ell}ext_{\ell}(T_{\mu(p)}(a))(x)=ext_{\ell}(\bigvee_{p\in\ell}T_{\mu(p)})(x)$
Also $ext_{\ell}$ is a subspace preserving map. Now we claim that $(id_X,ext_{\ell})$ is a continuous map in $\ell$-$\textbf{BSYM}$. To establish the claim it is necessary to show that $\models(x,ext_{\ell}(a))=ext_{\ell}(a)(x)=\models(x,a)=\models(id_X(x),a)$. For a given arrow $(f,\phi):P(S,\alpha)\rightarrow (X, A, \models)$ there exists a map, we define $\tilde{f}=f:(S,\alpha)\rightarrow Ext_{\ell}(X,A,\models)$ in $\ell$-BS. We show that the triangle in the Fig \ref{fig:coadjoint0} commute i.e., $(f,\phi)=(id_X, ext_{\ell})\circ (f, f^{-1})$. Here we observe that $id_X\circ f=f$. We have to prove that $f^{-1}\circ ext_{\ell}=\phi$.\\
Now 
$f^{-1}\circ ext_{\ell}(a)(x)= ext_{\ell}(a)(f(x))= \models(f(x), a)$, since $(f, \phi)$ is continuous we have $\models(f(x), a)= \models(x, \phi(a))$. 
Again $\models(x, \phi(a))=\phi(a)(x)$.

Hence $\xi_{(X,A,\models)}=(id_X, ext_{\ell})$ is the counit and as a result $Ext_{\ell}$ is the coadjoint to the functor $P$. 

\end{proof}
\end{thm}

\begin{figure}
 
\begin{center}
\includegraphics[width=100mm]{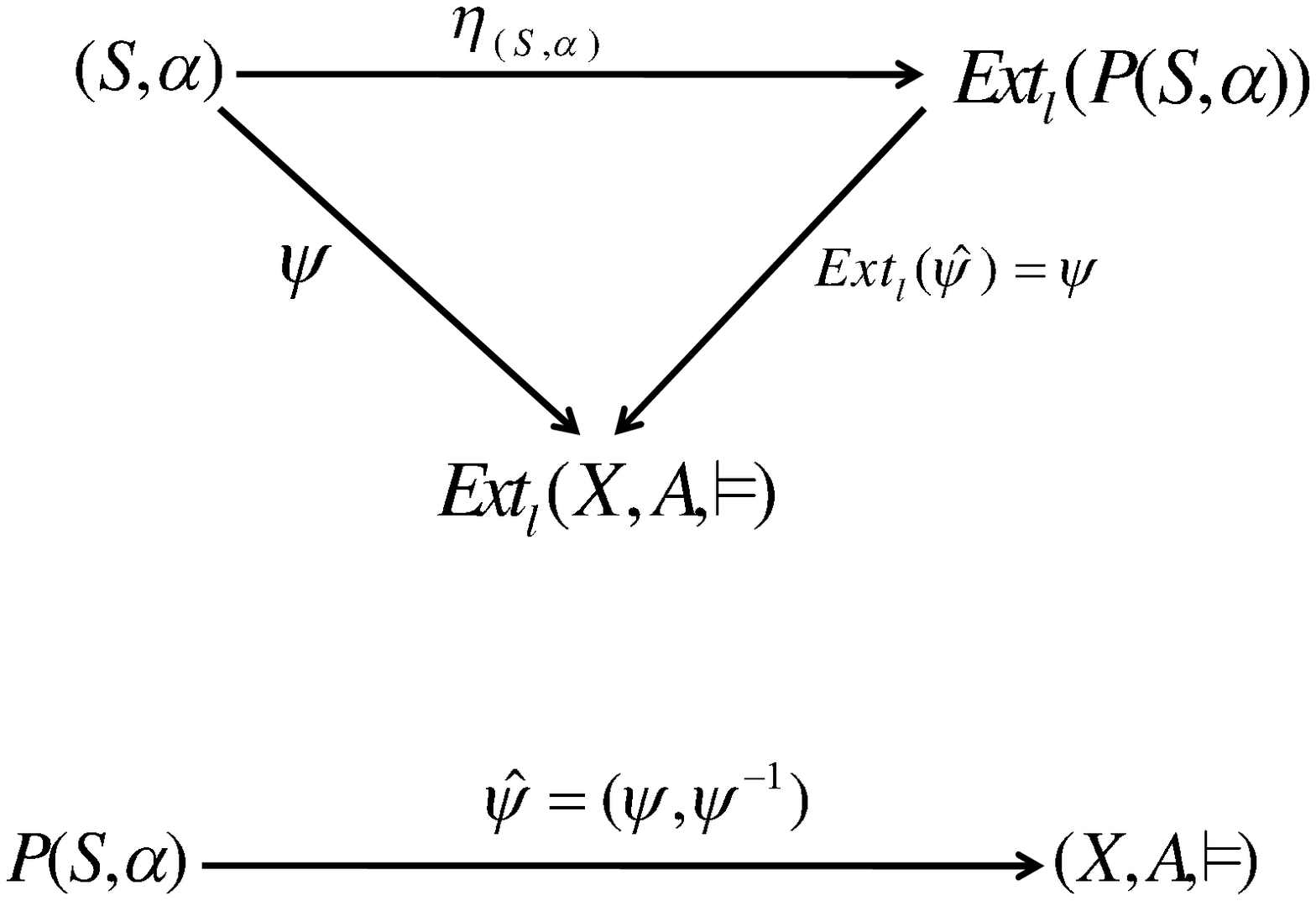}
\caption{Diagram of Unit}
\label{fig:adjoint0}
\end{center}
\end{figure}

Also $P$ is the adjoint to the functor $Ext_{\ell}$. Diagram of unit is shown in Fig. \ref{fig:adjoint0}.

\begin{thm}
\label{a1}
$Q$ is the adjoint to the functor $R$.

\begin{proof}

We prove the theorem by unit of the adjunction. Fig. \ref{fig:adjoint1} shows the diagram of unit.

\begin{figure}
 
\begin{center}
\includegraphics[width=100mm]{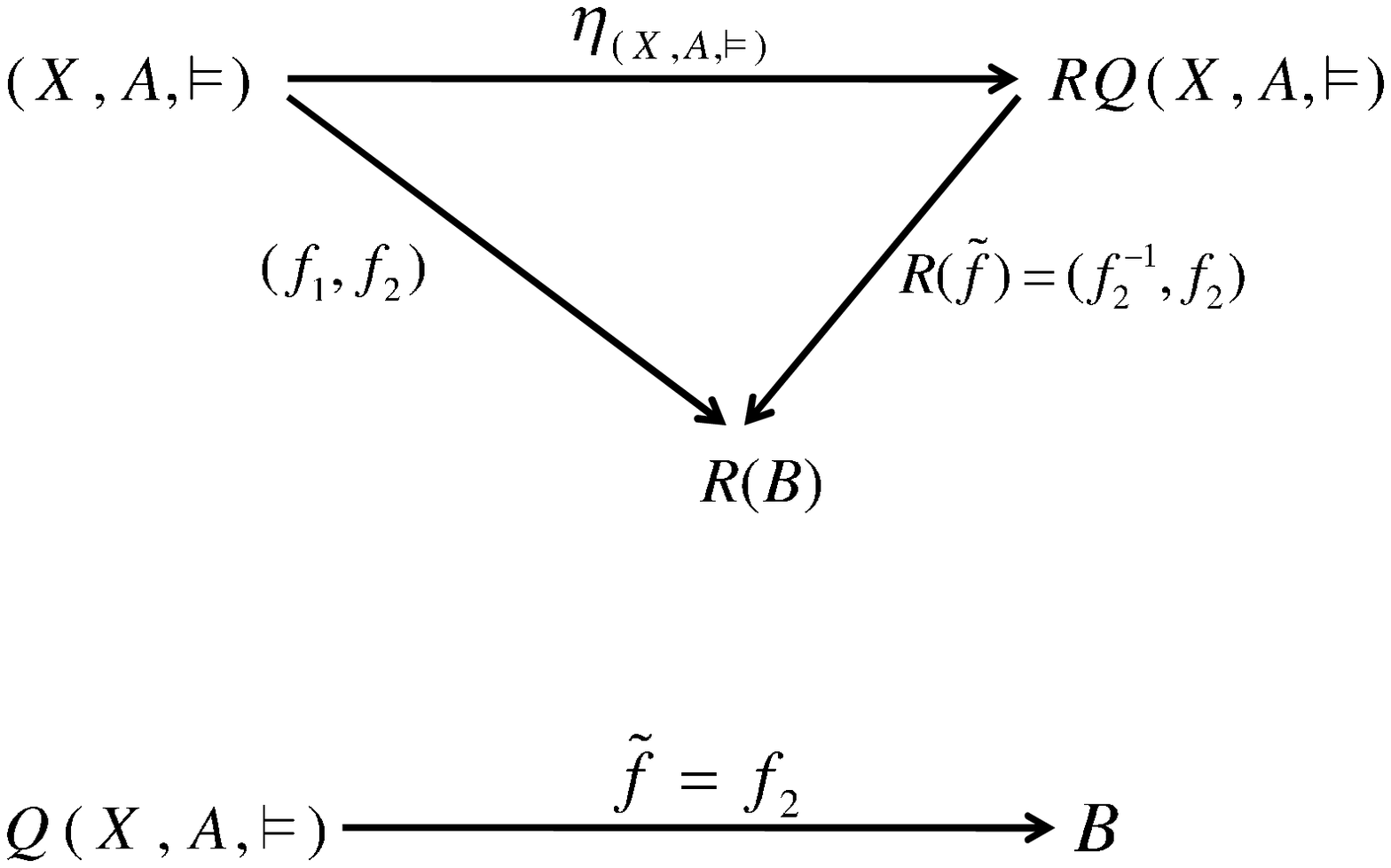}
\caption{Diagram of Unit}
\label{fig:adjoint1} 
\end{center}
\end{figure}
Here we recall that $R(B)=(Spec_{\ell}(B),B,\models_{(Spec_{\ell}(B)\times B)})$, where $\models_{(Spec_{\ell}(A)\times A)}(v,b)=v(b)$. So $RQ(X,A,\models)=R(A)=(Spec_{\ell}(A),A,\models_{(Spec_{\ell}(A)\times A)})$. Unit is defined by $\eta_{(X,A,\models)}=(f,id):(X,A,\models)\rightarrow RQ(X,A,\models)$ \\where
$f:X\rightarrow Spec_{\ell}(A)$ such that $f(x)=f_x$, $f_x:A\rightarrow \ell$. Where $f_x(a)=\models(x,a)$. Claim: for each $x\in X$, $f_x$ is an $\ell$-VL-algebra homomorphism. The claim is straightforward to check. The unit $\eta_{(X,A,\models)}=(f,id)$ is a continuous map of $\ell$-\textbf{BSYM}, since $\models(f(x),a)=f_x(a)=\models(x,a)=\models(id(x),a)$. For a given arrow $(f_1,f_2):(X,A,\models)\rightarrow R(B)$, we define $\tilde{f}=f_2$($f_2$ is an $\ell$-VL-algebra homomorphism). Now we shall show that the triangle in Fig. \ref{fig:adjoint1} commute i.e., $(f_1,f_2)=R(\tilde{f})\circ \eta_{(X,A,\models)}=(f_2^{-1},f_2)\circ (f,id)=(f_2^{-1}\circ f,id\circ f_2)$. It clearly shows that $id\circ f_2=f_2$. Now we show that $f_1=f_2^{-1}\circ f$. For each $x\in X$, $f_1(x)=f_2^{-1}\circ f(x)=f_x\circ f_2$ and for all $b\in B$, $f_x\circ f_2(b)=f_x(f_2(b))=\models(x,f_2(b))=\models(f_1(x),b)=f_1(x)(b)$. Therefore $f_x\circ f_2=f_1(x)$. Hence $f_2^{-1}\circ f=f_1$.

\end{proof}
\end{thm}

\begin{figure}
\begin{center}

\includegraphics[width=100mm]{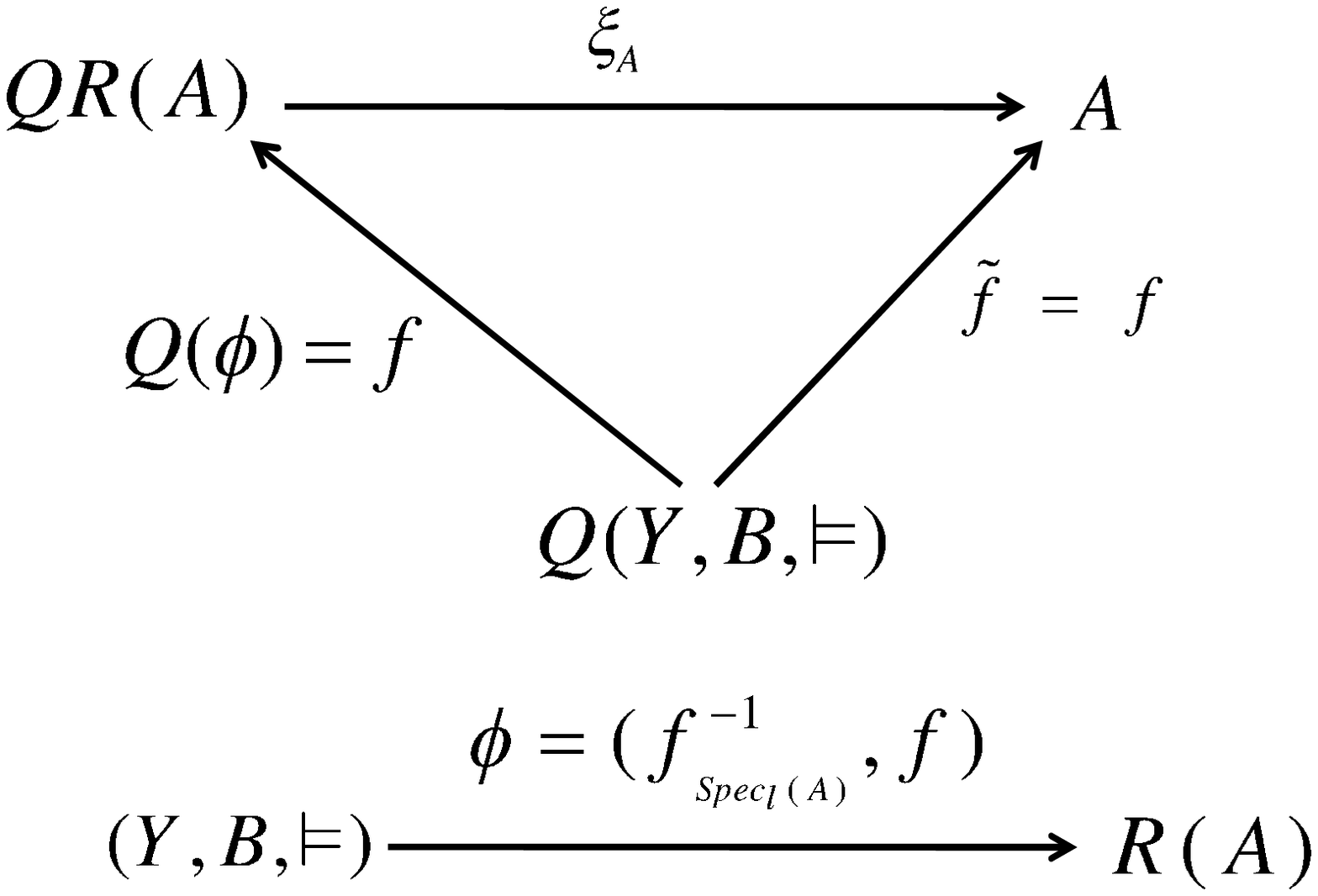}
\caption{Diagram of Counit}
\label{fig:coadjoint1}
\end{center}
\end{figure}

We can also prove $R$ is the coadjoint to the functor $Q$. Diagram of counit is shown in Fig. \ref{fig:coadjoint1}. 

For a given arrow $\tilde{f}$ in ($\ell$-VA)$^{op}$, we define $\phi=(f_{Spec_{\ell}(A)}^{-1},f)$ in $\ell$-\textbf{BSYM}. Where $f_{Spec_{\ell}(A)}^{-1}:Y\rightarrow Spec_{\ell}(A)$ defined by $f_{Spec_{\ell}(A)}^{-1}(y)=(v\circ f)(y)$, $v\in Spec_{\ell}(A)$ and $f$ is an $\ell$-VL-algebra homomorphism from $B$ to $A$ in ($\ell$-VA)$^{op}$.
It is easy to prove the triangle in Fig. \ref{fig:coadjoint1} commute i.e., $\xi_A\circ Q(\phi)=\tilde{f}$.
\begin{thm}
\label{01}
The category $\ell$-\textbf{BS} and $\ell$-\textbf{BSYM} are equivalent.
\begin{proof}
Let $id_{(X,A,\models)}$ and $id$ denote respectively the identity functors on $\ell$-\textbf{BSYM} and $\ell$-BS. We get two natural transformations $\xi$ and $\eta$ such that $\xi:P\circ ext_{\ell}\rightarrow id_{(X,A,\models)}$ and $\eta:id\rightarrow ext_{\ell}\circ P$. We will show $\xi_{(X,A,\models)}: P(ext_{\ell}(X,A,\models))\rightarrow (X,A,\models)$ is a natural isomorphism between two $\ell$-\textbf{BSYM}.
 Here we recall that $P(ext_{\ell}(X,A,\models))=(X,Cont(X,ext_{\ell}(A),\alpha)$ and $\xi_{(X,A,\models)}=(id, ext_{\ell})$. We have to show that $\xi_{(X, A,\models)}$  is a homeomorphism. Here $ext_{\ell}:A\rightarrow Cont(X,ext_{\ell}(A),\alpha)$ is an $\ell$-VL-algebra homomorphism. $id$ is of course injective and surjective.
 Only part we have to show $ext_{\ell}$ is an isomorphism. Let $ext_{\ell}(a_1)=ext_{\ell}(a_2)$. Then $ext_{\ell}(a_1)(x)=ext_{\ell}(a_2)(x) \Rightarrow \models(x, a_1)=\models(x, a_2)$. 
 So then for $r\in \ell$, $T_r(\models(x, a_1))=T_r(\models(x, a_2))\Rightarrow \models(x,T_r(a))=\models(x,T_r(a_2))$. Hence $ext_{\ell}(T_r(a_1))=ext_{\ell}(T_r(a_2))$. 
By proposition \ref{propo_1}(2) we get $T_r(a_1)\in P\Leftrightarrow T_r(a_2)\in P$. Claim: $T_r(a_1)=T_r(a_2)$. 
 If not, with out loss of generality we assume $T_r(a_1)\geq T_r(a_2)$. Let $\tilde{F}=\{p\in A:T_r(a_1)\leq p\}$. 
 Then $\tilde{F}$ is an $\ell$-filtrer. If we take $Y$ be the set of all $\ell$-filters $F'$ with $T_r(a_1)\in F'$ and $T_r(a_2) \not\in F'$. 
 Then by proposition \ref{propo_1}(1) there exists a prime $\ell$-filter $P'$ such that $T_r(a_1)\in P'$ but $T_r(a_2)\not\in P'$. 
 But this is a contradiction. Hence for all $r\in \ell$, $T_r(a_1)=T_r(a_2)$. 
 We know by definiotion of $\ell$-\textbf{VL}-algebra that $\displaystyle\bigwedge_{r\in \ell}(T_r(a_1)\leftrightarrow T_r(a_2))\leq a_1\leftrightarrow a_2$. But $\displaystyle\bigwedge_{r\in \ell}(T_r(a_1)\leftrightarrow T_r(a_2))=1$. 
 So $a_1\leftrightarrow a_2=1$. Therefore $a_1=a_2$ and hence $ext_{\ell}$ is injective.  Clearly $ext_{\ell}$ is surjective. Finally $\xi_{(X,A,\models)}$ satifies the continuity condition i.e., $\models(x,ext_{\ell}(a))=ext_{\ell}(a)(x)=\models(x,a)=\models(id(x),a)$. 
Therefore $\xi_{(X,A,\models)}$ is an isomorphism. As a result $\xi$ is a natural isomorphism.\\
Now we will show that $\eta$ is a natural isomorphism. Here we recall that $Ext_{\ell}(P(S,\alpha))=(S,ext_{\ell}(Cont(S,\alpha),\alpha')$. 
We define $\eta_{(S,\alpha)}:(S,\alpha)\rightarrow Ext_{\ell}\circ P(S,\alpha)$ by $\eta_{(S,\alpha)}(s)(f)=f(s)$, $f\in Cont(S,\alpha)$. 
Where the map $\alpha':Subalg(\ell)\rightarrow Subsp(S,ext_{\ell}(Cont(S,\alpha))$ such that $\alpha'(m)=(S,ext_{m}(S,Cont(S,\alpha))$. 
Now if $x\in \alpha(m)$, for $m\in Subalg(\ell)$ then $\eta_{(S,\alpha)}(x)(f)=f(x)\in \alpha'(m)$. It can be shown that $\eta$ is a natural isomorphism.

\end{proof}
\end{thm}
\begin{thm}
\label{02}
$\ell$-VA is dually equivalent to $\ell$-\textbf{BSYM}.
\begin{proof}
We have two natural transformations $\xi$ and $\eta$ such that $\xi_A=id:Q(R(A))\rightarrow A$ and $\eta_{(X,A,\models)}:(X,A,\models)\rightarrow (Spec_{\ell}(A),A,\models_{(Spec_{\ell}(A)\times A)})$. It is clear that $\xi$ is natural isomorphism. We shal show that $\eta_{(X,A,\models)}$ is a natural isomorphism between two boolean systems. Here $\eta_{(X,A,\models)}=(f,id_{A})$ such that \\
\begin{itemize}
\item $f:X\rightarrow Spec_{\ell}(A)$ is a set map.
\item $id_{A}: A\rightarrow A$ is an $\ell$-VL-algebra homomorphism.
\end{itemize}
We have to show $\eta_{(X,A,\models)}$ is a homeomorphism. First we show that $f$ is bijective.
Claim: $f$ is injective and surjective. Let $x_1\neq x_2$ in $X$.Then by Definition \ref{bsym} we have $\models(x_1,a)\neq \models(x_2,a)$, for some $a\in A$. Therefore $f(x_1)(a)\neq f(x_2)(a)$,for some $a\in A$. As a result $f$ is injective. From the definition of $f$, we can say that $f$ is surjective also. Hence prove our claim.\\
Finally we observe that $f(x)(a)=f_x(a)=\models(x,a)$ and $f_x(a)=\models_{(Spec_{\ell}(A)\times A)}(f(x),a)$. Therefore $\models_{(Spec(A)\times A)}(f(x),a)=\models(x,id_{A}(a))$. Hence $\eta_{(X,A,\models)}$ is an isomorphism and therefore $\ell$-VA is dually equivallent to $\ell$-\textbf{BSYM}.
\end{proof}
\end{thm}

\begin{thm}
\label{03}
$\ell$-VA is dually equivallent to $\ell$-\textbf{BS}.
\begin{proof}
As adjunctions can be composed, hence composition of equivalences of Theorems \ref{01} and \ref{02} shows the result.
\end{proof}
\end{thm}
\section{ $\ell$-\textbf{ML}-algebras, $\ell$-relational systems, $\ell$-relational space and their Categorical interconnections}
\label{sec:3}
\textbf{$\ell$-\textbf{ML}-algebras}\\
We now introduce the notion of $\ell$-\textbf{ML}-algebras(for details see \cite{maruyama2009algebraic}). 
\begin{defn}
An $\ell$-\textbf{ML}-algebra is an algebraic structure $(A, \wedge, \vee, \rightarrow, T_r(r\in \ell), \Box, 0, 1)$ such that the following hold.
\begin{enumerate}[(i)]
\item $(A, \wedge, \vee, \rightarrow, T_r(r\in \ell), 0,1)$ is an $\ell$-\textbf{VL}-algebra.
\item $\Box(a_1\wedge a_2)=\Box a_1\wedge\Box a_2$ and $\Box 1=1$
\item for all $r\in \ell$, $U_r(\Box a)=\Box U_r(a)$. Where $U_r(a)=\bigvee\{T_{r_1}(a)| r\leq r_1\}$.
\end{enumerate}
\end{defn}
\begin{defn}(\cite{maruyama2011dualities})
$\ell$-\textbf{ML}-algebra homomorphism is defined as a homomorphism of $\ell$-\textbf{VL}-algebras and satify the operation $\Box$.
\end{defn}
For an $\ell$-\textbf{ML}-algebra $A$, $Spec_{\ell}(A)$ is same as given in Definition \ref{spec}. 
\begin{defn}(\cite{maruyama2011dualities})\label{Spec_R}
A binary relation $R$ on $Spec_{\ell}(A)$ is defined as follows: $fRg \Leftrightarrow \forall r\in \ell, \forall a\in A, f(\Box a)\geq r \Rightarrow g(a)\geq r$. Then $(Spec_{\ell}(A), R, e)$ is an $\ell$-valued canonical model of $A$, where $e$ is a Kripke $\ell$-valuation on $(Spec_{\ell}(A), R)$  such that  $e(f,a)=f(a)$, $\forall f\in Spec_{\ell}(A)$.

\end{defn}
\begin{prop}(\cite{maruyama2011dualities})
\label{Spec_P}
The $\ell$-valued canonical model $(Spec_{\ell}(A), R, e)$ of $A$ is an $\ell$-valued Kripke model. In other words, $e(f,\Box a)=f(\Box a)=\bigwedge\{g(a) | fRg\}$.
\end{prop}
\begin{prop}(\cite{maruyama2011dualities})
For an $\ell$-\textbf{ML}-algebra $A$, the Boolean algebra $\mathcal{B}(A)$ is a modal algebra.
\end{prop}
\textbf{$\ell$-relational systems.}
\begin{defn}
\label{lrsym}
An $\ell$-relational systems is a triple $(X,A,\models)$ where $X$ is a nonempty set, $A$ is an $\ell$-\textbf{ML}-algebra and $\models$ is an $\ell$-satisfaction relation from $X$ to $A$ such that the following hold.
\begin{enumerate}
\item $\models$ satisfies both the join and finite meet interchange law.
\item $\models(w,\Box a)=\wedge\{\models(w',a)|wR_0w'\}$, $R_0$ is a binary relation on $X$.
\item $\models(w,T_r(a))=T_r(\models(w,a))$.
\item $\models(x,\bot)=\bot$, $\models(x,\top)=\top$.
\item $\models(x,a\rightarrow b)=\models(x,a)\rightarrow \models(x,b)$.
\end{enumerate}
\end{defn}

\subsection{Categories and their functorial relationships}
Following the Definition \ref{cbsym} we define a category $\ell$-\textbf{RSYM} as objects are $\ell$-relational systems and arrows are continuous functions.
\begin{defn}(\cite{maruyama2011dualities})
\label{lrs}
We use the definition of the category $\ell$-\textbf{RS} given in \cite{maruyama2009algebraic} and it was defined as follows.
\begin{enumerate}
\item Objects: $(S,\psi,R)$ where $(S,\psi)$ is an object in $\ell$-BS and $R$ is binary relation on $S$ which has the following properities.
\begin{itemize}
\item if $\forall f \in Cont(S,\psi)$, $(\Box_Rf)(x)=1\Rightarrow f(y)=1$ then $(x,y)\in R$.
\item if $X$ is a clopen subset of $S$ then $R^{-1}[X]$ is a clopen subset of $S$.
\item for any $m\in Subalg(\ell)$, if $x\in \psi(m)$ then $R[x]\subset \psi(m)$.
\end{itemize}
\item Arrows: $f:(S_1,\psi_1,R_1)\rightarrow (S_2,\psi_2,R_2)$ in $\ell$-RS is an arrow(morphism) $f:(S_1,\psi_1)\rightarrow (S_2,\psi_2)$ in $\ell$-BS having the following properities.
\begin{itemize}
\item if $sR_1t$ then $f(s)R_2f(t)$.
\item if $f(s_1)R_2s_2$ then there exists $t_1\in S_1$ such that $s_1R_1t_1$ and $f(t_1)=s_2$.
\end{itemize}
\end{enumerate}
\end{defn}
\begin{defn}
\label{rel}
For an $\ell$-\textbf{ML}-algebra $A$, We define a binary relation $R$ on $(X,ext_{\ell}(A))$ as follows.\\
$xRy \Leftrightarrow \forall r\in \ell, \forall a\in A,  ext_{\ell}(\Box a)(x)\geq r \Rightarrow ext_{\ell}(a)(y)\geq r$.
\end{defn}

\begin{figure}

\begin{center}
\includegraphics[width=100mm]{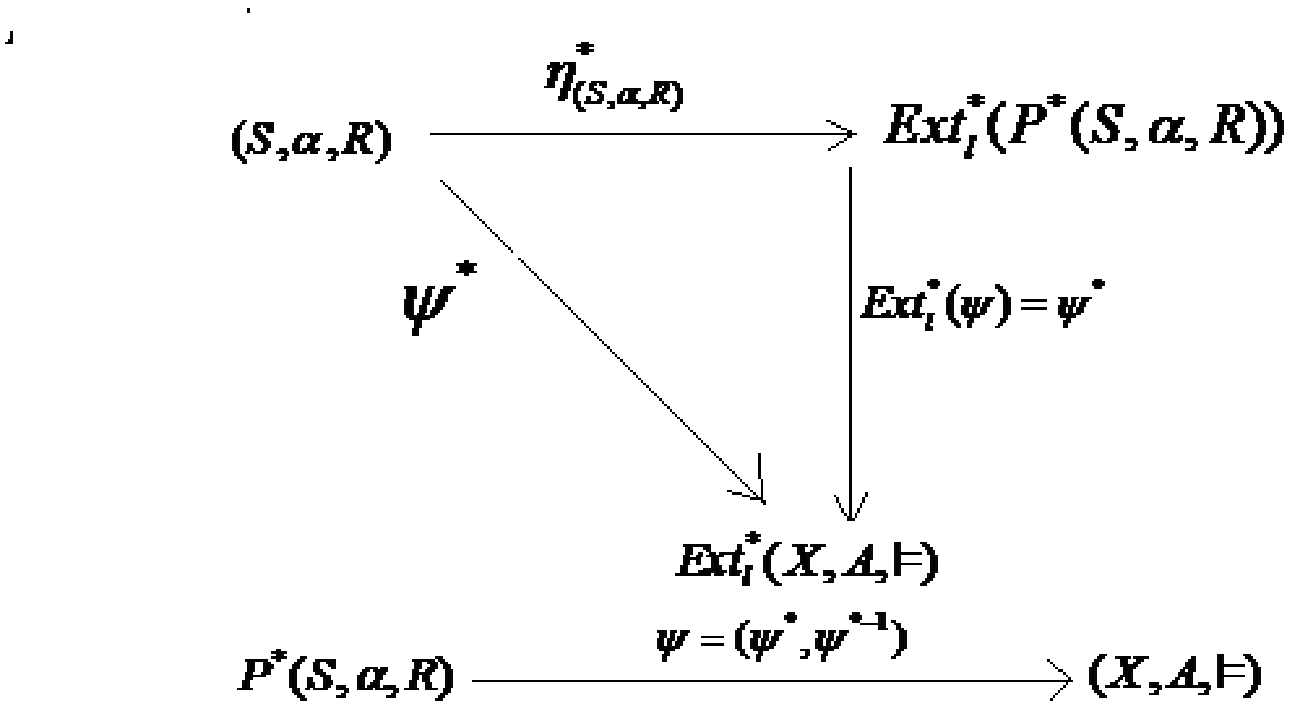}
\caption{Diagram of Unit}
\label{fig:adjoint_11}
\end{center}
\end{figure}

\begin{defn}
$Ext^*_{\ell}$ is a functor from $\ell$-\textbf{RSYM} to $\ell$-\textbf{RS} defined as follows.
\begin{enumerate}
\item acts on an object $(X,A,\models_1)$ as $Ext^*_{\ell}(X,A,\models_1)=((X,ext_{\ell}(A)),\alpha,R_{\Box})$,
\item for an arrow $(f,\phi):(X,A,\models_1)\rightarrow (Y,B,\models_2)$ in $\ell$-\textbf{RSYM} define $Ext^*_{\ell}(f,\phi)=f$.
\end{enumerate}
\end{defn}
Now we verify the well-definedness of $Ext^*_{\ell}$.
\begin{lem}
$(X,ext_{\ell}(A),\alpha,R_{\Box})$ is an object of $\ell$-\textbf{RS}.
\begin{proof}
We show that, if for all $f\in Cont(X,ext_{\ell}(A),\alpha)(\Box_Rf)(x)=1\Rightarrow f(y)=1$, then $xR_{\Box}y$.
We prove the contrapositive statement. Suppose $(x,y)\notin R_{\Box}$. Then there exists $r\in \ell$ and $a\in A$ such that $ext_{\ell}(\Box a)(x)\geq r \Rightarrow ext_{\ell}(a)(y)\ngeq r$. Now $U_rext_{\ell}(\Box a)(x)=1 \Rightarrow ext_{\ell}(U_r(\Box a))(x)=1$. But $ext_{\ell}(U_r(a))(y)\neq 1$. Define $f:(X,ext_{\ell}(A),\alpha)\rightarrow (\ell,\alpha_{\ell})$ by $f(x)=ext_{\ell}(U_r(a))(x)$. We have $(\Box_Rf)(x)=\wedge \{f(y):xR_{\Box} y\}=\wedge \{ext_{\ell}(U_r(a)(y):xR_{\Box}y\}=ext_{\ell}(\Box U_r(a))(x)=1$. But $f(y)=ext_{\ell}(U_r(a))(y)\neq 1$. We know $ext_{\ell}(a)\in Cont((X,ext_{\ell}(A)),\alpha)$. So by definition of $f$, we have $f\in Cont((X,ext_{\ell}(A)),\alpha)$.\\
Now we verify the second position in the object section of Definition \ref{lrs}.  For $r\in \ell$, $(ext_{\ell}(a))^{-1}\{r\}=(T_r\circ ext_{\ell}(a))^{-1}\{1\}$ is a clopen set i.e., both open and closed (since $T_r\circ ext_{\ell}(a)\in Cont(X,ext_{\ell}(A),\alpha)$). Now $R_{\Box}^{-1}[(ext_{\ell}(a)^{-1}\{r\}]=R_{\Box}^{-1}[(T_r\circ ext_{\ell}(a))^{-1}\{0\}]$, is clopen in $X$.\\
After this, we verify the third position in the object section of Definition \ref{lrs}. Here $\alpha(s)=(X,ext_s(A))$, where $s$ is a subalgebra of $\ell$. Let $z\in (X,ext_s(a))$ and $R[z]\setminus \alpha(s)\neq \phi$. Then for $w\in R[z]\setminus \alpha(s)$ we have $ext_s(a)(w)\notin s$.Define $ext_s(a)(w)=r$. Now for $w'\in (X,ext_{\ell}(A))$

\begin{equation*}
ext_{\ell}(T_r(a)\rightarrow a)(w')=\begin{cases}
                                                      1, &\text{if}\ ext_{\ell}(a)(w')\neq r  \\
                                                      r, &\text{if}\ ext_{\ell}(a)(w')=r
\end{cases}
\end{equation*}                   
Now $ext_{s}(\Box(T_r(a)\rightarrow a))(z)=\models(z,\Box(T_r(a)\rightarrow a))=\wedge \{\models(z',T_r(a)\rightarrow a)|zR_{\Box}z'\}=\wedge \{ext_s(T_r(a)\rightarrow a)(z')|zR_{\Box}z'\}=r$. But this contradicts our assumption that $ext_{s}(\Box(T_r(a)\rightarrow a))(z)\in s$ . Therefore if $z\in \alpha(s)$ then $R[z]\subset\alpha(s)$. 

\end{proof}
\end{lem}
\begin{lem}
For an arrow $(f,\phi):(X,A,\models_1)\rightarrow (Y,B,\models_2)$ in $\ell$-\textbf{RSYM}, $Ext^*_{\ell}(f,\phi)$ is an arrow in $\ell$-\textbf{RS}.
\begin{proof}
We verify $Ext^*_{\ell}(f,\phi)=f:(X,ext_{\ell}(A)\alpha_1,R_{\Box_1})\rightarrow (Y,ext_{\ell}(B),\alpha_2,R_{\Box_2})$ is an arrow in $\ell$-\textbf{RS}.
Here we note that $f$ is an arrow in $\ell$-\textbf{BS}. Assume $xR_{\Box_1}y$. Claim $f(x)R_{\Box_2}f(y)$. By Definition \ref{rel} we have $ext_{\ell}(\Box a)(x)\geq r\Rightarrow ext_{\ell}(a)(y)\geq r$. Now if for all $b\in B$ and $r_1\in \ell$, $ext_{\ell}(\Box b)f(x)\geq r_1$, then $\models_2(f(x),\Box b)\geq r_1$. By Definition \ref{lrsym} we have $\wedge\{\models_2(f(y),b)|xR_{\Box_1}y\}\geq r$. 
This shows that $ext_{\ell}(b)(f(y))\geq r_1$.\\
We next verify that $Ext^*_{\ell}(f,\phi)$ satisfies the condition 2 in the arrow section of Definition \ref{lrs}. Assume $f(x_1)R_{\Box_2}x_2$. Define $Ext_2(f^*,\phi^*):(X,ext_2(\mathcal{B}(A_1)),\alpha_1^*,R_{\Box_1}^*)\rightarrow (Y,ext_2(\mathcal{B}(A_2)),\alpha_2^*,R_{\Box_2}^*)$ by $Ext_2(f^*,\phi^*)=f^*$, where $f^*(x)=f(x)$ for $x\in (X,ext_2(\mathcal{B}(A_1)))$.
It can be shown that $Ext_2(f^*,\phi^*)$ is an arrow in 2-\textbf{RS}. We have if $f^*(x_1)R_{\Box_2}^*x_2$ then there is $y_1$ in $(X,ext_2(\mathcal{B}(A_1)))$ such that $x_1R_{\Box_1}^*y_1$ and $f^*(y_1)=x_2$.
Now $ext_{\ell}(a_1)(y_1)=r\Leftrightarrow ext_2(T_r(a_1))(y_1)=1$. We claim $x_1R_{\Box_1}y_1$ and $f(y_1)=x_2$. If $ext_{\ell}(\Box a)(x_1)\geq r$ then $T_1\circ (ext_{\ell}(\Box U_r(a)))(x_1)=1$.
Therefore $ext_{\ell}(\Box T_1(U_ra))(x_1)=1$. Since $x_1R_{\Box_1}^*y_1$, we have $ext_{\ell}(U_ra)(y_1)=1\Rightarrow ext_{\ell}(a)(y_1)\geq r$. Therefore we have $x_1R_{\Box_1}y_1$.
Let $ext_{\ell}(b)(f(y_1))=r$. Then $ext_2(T_r(b))(f^*(y_1))=1\Rightarrow ext_2(T_r(b))(x_2)=1$. Hence $ext_{\ell}(b)(x_2)=r$. Therefore 
\begin{equation*}
ext_{\ell}(b)(x_2)=ext_{\ell}(b)(f(y_1))\\
\Rightarrow \models_2(x_2,b)=\models_2(f(y_1),b)\\
\Rightarrow f(y_1)=x_2
\end{equation*}

\end{proof}
\end{lem}
\begin{defn}
$P^*$ is a functor from $\ell$-\textbf{RS} to $\ell$-\textbf{RSYM} defined as follows. \\
\begin{enumerate}[(i)]
\item $P^*$ acts on an object $(S,\alpha,R_{\Box})$ in $\ell$-\textbf{RS} as $P^*(S,\alpha,R_{\Box})=(S,(Cont(S,\alpha),\Box_R),\models)$ where $\models(s,v)=v(s)$. $(Cont(S,\alpha),\Box_R)$ is an $\ell$-\textbf{ML}-algebra(for details see( \cite{maruyama2011dualities})).
\item $P^*$ acts on an arrow(morphism) $f:(S_1,\alpha_1,R_1)\rightarrow (S_2,\alpha_2,R_2)$ as $P^*(f)=(f,f^{-1})$.
\end{enumerate}
 Where 
\begin{itemize}
\item $f:S_1\rightarrow S_2$ a set function.
\item $f^{-1}:(Cont(S_2,\alpha_2),\Box_{R_1})\rightarrow (Cont(S_1,\alpha_1),\Box_{R_2})$ is an $\ell$-\textbf{ML}-algebra homomorphism. Where $f^{-1}(v)=v\circ f$, $v\in Cont(S_2,\alpha_2)$.
\end{itemize}
\end{defn}
\begin{defn}
$Q^*$ is a functor from $\ell$-\textbf{RSYM} to ($\ell$-MA)$^{op}$ defined as follows.
\begin{enumerate}
\item $Q^*$ acts on an object $(X,A,\models)$ in $\ell$-\textbf{RSYM} as $Q^*(X,A,\models)=A$.
\item $Q^*$ acts on an arrow(morphism) $(f,\phi):(X,A,\models_1)\rightarrow (Y,B,\models_2)$ in $\ell$-\textbf{RSYM} as $Q^*(f,\phi)=\phi$
\end{enumerate}
\end{defn}
\begin{lem}
For an $\ell$-\textbf{ML}-algebra $A$, $(Spec_{\ell}(A),A,\models_{(Spec_{\ell}(A)\times A)})$ is an object in $\ell$-\textbf{RSYM}.
\begin{proof}
Here we use $\models$ rather than $\models_{(Spec_{\ell}(A)\times A)}$. We define $\models(v,a)=v(a)$, $v\in Spec_{\ell}(A)$. Clearly $Spec_{\ell}(A)$ is a set and $A$ is a frame. Next we show that\\
\begin{enumerate}[(i)]
\item $\displaystyle\models(v,\bigvee_{\gamma\in \Gamma}a_{\gamma})=v(\bigvee_{\gamma\in\Gamma}a_{\gamma})=\bigvee_{\gamma\in\Gamma}v(a_{\gamma})=\bigvee_{\gamma\in\Gamma}\models(v,a_{\gamma})$. Also $\models(v,a_1\wedge a_2)=v(a_1)\wedge v(a_2)$.
\item 
\begin{align*}
\models(v,\Box a)&=v(\Box a)\\
                          &=\wedge\{u(a): vR_{\Box}u\} \text{( using the Proposition \ref{Spec_P})}\\
                          &=\wedge\{\models(u,a) | vR_{\Box}u\}
\end{align*}
\item $\models(v,T_r(a))=v(T_r(a))=T_r(v(a))=T_r(\models(v,a))$.
\item $\models(v,a\rightarrow b)=v(a\rightarrow b)=v(a)\rightarrow v(b)=\models(v,a)\rightarrow \models(v,b)$
\end{enumerate}
Therefore $(Spec_{\ell}(A),A,\models_{(Spec_{\ell}(A)\times A)})$ is an $\ell$-relational system.
\end{proof}
\end{lem}
It is easy to prove that $(\psi^{-1},\psi)$ is continuous whenever $\psi$ is an $\ell$-\textbf{ML}-algebra homomorphism.
\begin{defn}
$R^*$ is a functor from ($\ell$-MA)$^{op}$ to $\ell$-\textbf{RSYM} defined thus.
\begin{enumerate}
\item $R^*$ acts on an object $A$ in $\ell$-MA as $R^*(A)=(Spec_{\ell}(A),A,\models_{(Spec_{\ell}(A)\times A)})$.
\item $R^*$ acts on an arrow(morphism) $\psi:A\rightarrow B$ in ($\ell$-MA)$^{op}$ as $R^*(\psi)=(\psi^{-1},\psi)$.
\end{enumerate}
\end{defn}
By the above lemma it can be shown that $R^*$ is indeed a functor.
\begin{thm}
$Ext^*_{\ell}$ is the co-adjoint to the functor $P^*$.
\begin{proof}
We will prove the theorem showing counit of the adjunction. Fig \ref{fig:coad22} shows the diagram of counit.
\begin{figure}

\begin{center}
\includegraphics[width=100mm]{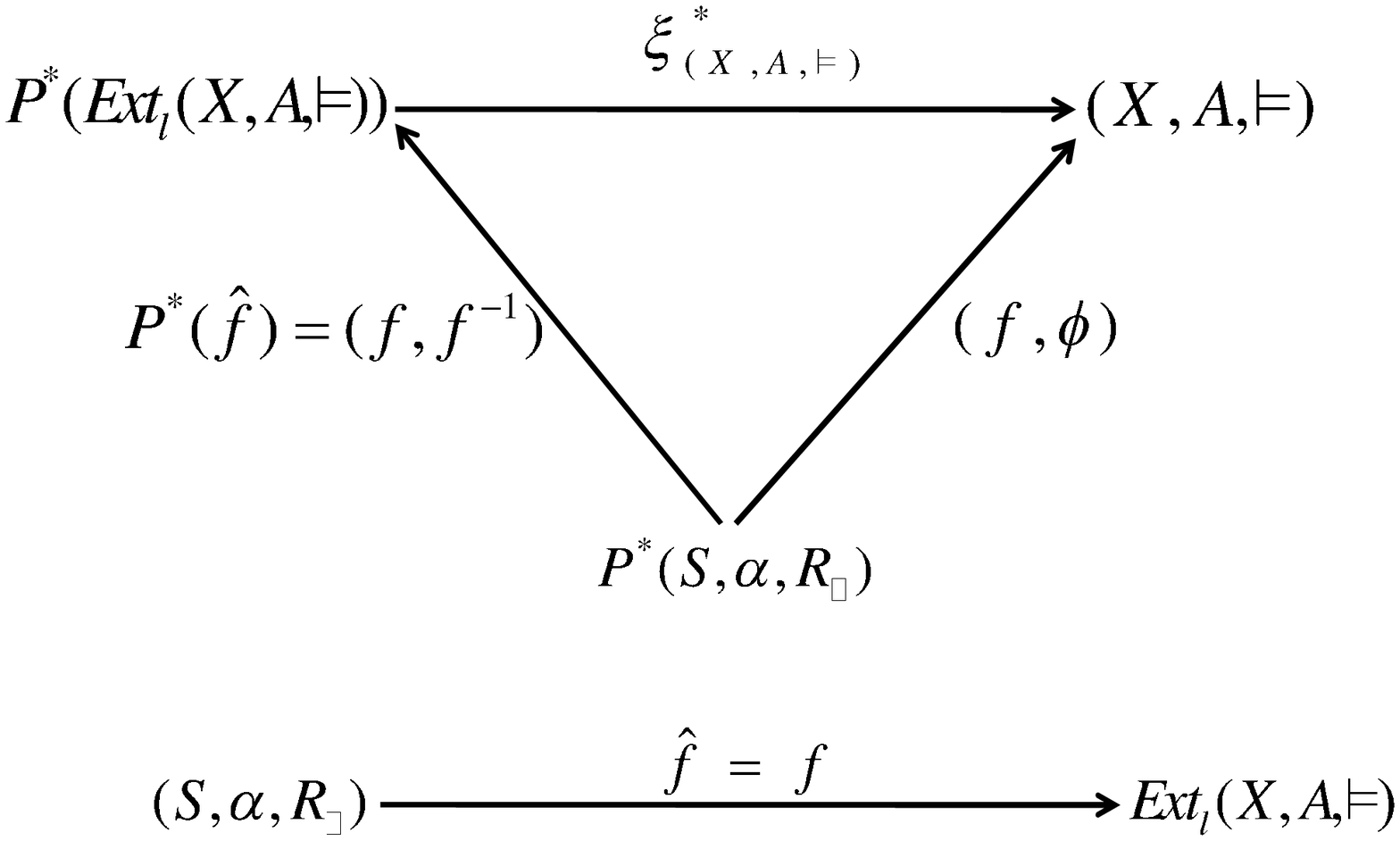}
\caption{Diagram of Counit}
\label{fig:coad22}
\end{center}
\end{figure}
Recall that $P^*(S,\alpha,R_{\Box})=(S,(Cont(S,\alpha),\Box_R),\models)$ and $Ext_{\ell}(X,A,\models)=(X,ext_{\ell}(A),\alpha,R_{\Box})$. So $P^*(Ext_{\ell}(X,A,\models))=(X,(Cont(X,ext_{\ell}(A),\alpha),\Box_R),\models)$. 

We now show that counit $\xi^*_{(X,A,\models)}=(id^*_X, ext^*_{\ell}) : P^*(Ext^*_{\ell}(X,A,\models))\rightarrow (X,A,\models)$ is a continuous map of $\ell$-relational system. Here 
\begin{enumerate}[(i)]
\item $id_X: X\rightarrow X$ is clearly a set function.
\item $ext^*_{\ell}: A\rightarrow (Cont(X,ext_{\ell}(A),\alpha),\Box_{R})$ is an $\ell$-\textbf{ML}-algebra homomorphism. Where $ext^*_{\ell}(a)=ext_{\ell}(a)$, $\forall a\in A$.
\end{enumerate}
We know that $ext^*_{\ell}$ is an $\ell$-VL-algebra homomorphism. So only part we have to show is that it preserves the operations $\Box$ i.e., $ext^*_{\ell}(\Box a)=\Box(ext^*_{\ell}(a))$. Now 
$ext^*_{\ell}(\Box a)(x)=\models(x,\Box a)=\wedge\{\models(y,a):xR_{\Box}y\}$. Again $\Box(ext^*_{\ell}(a))(x)=\wedge\{ext_{\ell}(a)(y): xR_{\Box}y\}$. So $ext^*_{\ell}$ is an $\ell$-\textbf{ML}-algebra homomorphism.
To prove the continuity of $\xi^*_{(X, A,\models)}$ it is enough to show that $\models(id^*_X(x),a)=\models(x,ext^*(a))$. We see 
$\models(id^*_X(x),a) =\models(x,a)=ext_{\ell}(a)(x)=\models(x,ext_{\ell}(a))=\models(x,ext^*(a))$.\\
Next we prove that the triangle in the Fig. \ref{fig:coad22} commute i.e., for a given arrow $(f,\varphi): P^*(S,\alpha,R_{\Box})\rightarrow (X,A,\models)$ there is an arrow, we define $\hat{f}=f$ such that 
$(f,\varphi)=\xi^*_{(X,A,\models)}\circ P^*(\hat{f})$.\\ Now 
\begin{align*}
(f,\varphi)&=(id^*_X, ext^*_{\ell})\circ (f, f^{-1})\\
              &=(id^*_X\circ f, f^{-1}\circ ext^*_{\ell})
\end{align*}
It is clear that $id^*_X\circ f=f$. Only part we have to show $\varphi=f^{-1}\circ ext^*_{\ell}$. Now as $(id^*_X, ext^*_{\ell})$ is continuous, so $\models(id^*_X(x), a)=\models(x, ext^*_{\ell}(a))$. Therefore $ext^*_{\ell}(a)=a$.\\
Now

\begin{align*}
f^{-1}\circ ext^*_{\ell}(a) &=f^{-1}(a)\\
&=\varphi(a) \text{( as} (f,\varphi)  \text{ is continuous)}
\end{align*}

 Hence $\xi^*_{(X, A, \models)}$ is the counit and as a result $Ext^*_{\ell}$ is the coadjoint to the functor $P^*$.
 
\end{proof}
\end{thm}
It can also be shown that $P^*$ is the adjoint to the functor $Ext^*_{\ell}$.
Fig. \ref{fig:adjoint_11} shows the diagram of unit.
For a given arrow $\psi^* : (S, \alpha, R)\rightarrow Ext^*_{\ell}(X, A, \models)$ there is an arrow, we define $\psi : P^*(S, \alpha, R)\rightarrow (X, A, \models)$ such that $Ext^*_{\ell}(\psi)=\psi^*$. It is easily check that the triangle in the Fig. \ref{fig:adjoint_11} commute i.e., $\eta^*_{(S,\alpha,R)}\circ Ext^*_{\ell}(\psi)=\psi^*$. 
\begin{thm}
\label{3}
The category $\ell$-\textbf{RSYM} is equivalent to the category $\ell$-\textbf{RS}.
\begin{proof}
Let $id^*_{(X,A,\models)}$ and $id^*$ denote respectively the identity functors on $\ell$-\textbf{RSYM} and $\ell$-\textbf{RS}. $\xi^*$ and $\eta^*$ are two natural transformations such that for an object $(X,A,\models)$ in $\ell$-\textbf{RSYM}, 
$\xi^*_{(X,A,\models)}: P^*(Ext^*_{\ell}(X, A, \models))\rightarrow (X, A, \models)$ and for an object $(S,\alpha,R_{\Box})$ in $\ell$-\textbf{RS}, $\eta^*_{(S,\alpha,R_{\Box})}: (S, \alpha, R_{\Box})\rightarrow Ext^*_{\ell}(P^*(S, \alpha, R_{\Box}))$.
 Now it is enough to show that $\xi^*$ and $\eta^*$
are natural isomorphism. $\xi^*_{(X,A,\models)}$ is a natural transformation between two $\ell$-\textbf{RSYM}. We here mention that $\xi^*_{(X,A,\models)}$ is almost same as $\xi_{(X,A,\models)}$ in the proof of Theorem \ref{01}.
 So by Theorem \ref{01}, $\xi^*_{(X,A,\models)}$ is an isomorphism and hence $\xi^*$ is a natural isomorphism.\\ 
Now it is left to prove that $\eta^*$ is a natural isomorphism. $\eta^*_{(S,\alpha,R_{\Box})}: (S,\alpha,R_{\Box})\rightarrow (S,ext^*_{\ell}(Cont(S,\alpha),\Box_R),\alpha',R'_{\Box})$ defined by $\eta^*_{(S,\alpha,R_{\Box})}(s)(f)=f(s)$, $f\in Cont(S,\alpha)$.
 Since $ext^*_{\ell}(f)(s)=f(s)$, so $\eta^*_{(S,\alpha,R_{\Box})}$ is well defined. 
Here $\eta^*_{(S,\alpha,R_{\Box})}$ is almost same as $\eta_{(S,\alpha)}$ in the proof of the Theorem \ref{01}.
 So by Theorem\ref{01} $\eta^*_{(S,\alpha,R_{\Box})}$ is an isomorphism in the category $\ell$-BS. We have to show that $\eta^*_{(S,\alpha,R_{\Box})}$ and $\eta^{*}{-1}_{(S,\alpha,R_{\Box})}$ satisfy the conditions 1 and 2 in the arrow part of Definition \ref{lrs}. 
 Assume for any $s_1, s_2\in S$, $s_1R_{\Box}s_2$. Then for any $r\in \ell$ and $f\in Cont(S,\alpha)$, $\eta^*_{(S,\alpha,R_{\Box})}(s_1)(\Box_R f)=ext^*_{\ell}(\Box_Rf)(s_1)\geq r\Rightarrow (\Box_Rf)(S_1)\geq r$.
 Now $(\Box_Rf)(s_1)=\wedge\{f(s_3):s_1R_{\Box}s_3\}$. Since $s_1R_{\Box}s_2$, we have $f(s_2)\geq r$. Therefore $ext_{\ell}(s_2)\geq r$ and hence $\eta^*_{(S,\alpha,R_{\Box})}(s_1)R'_{\Box}\eta^*_{(S,\alpha,R_{\Box})}(s_2)$.
 Again we observe that if $(s_1,s_2)\notin R_{\Box}$ then by object part of Definition \ref{lrs} there exists $f\in Cont(S,\alpha)$ such that $(\Box_{R_{\Box}}f)(s_1)=1$ but $f(s_2)\neq 1$. Therefore $ext^*_{\ell}(\Box_{R_{\Box}}f)(s_1)=1$ and $ext^*_{\ell}(f)(s_2)\neq 1$. 
 Therefore
  $(\eta^*_{(S,\alpha,R_{\Box})}(s_1), \eta^*_{(S,\alpha,R_{\Box})}(s_2))\notin R'_{\Box}$. So we get for any $s_1, s_2\in S$, $s_1R_{\Box}s_2$ iff $\eta^*_{(S,\alpha,R_{\Box})}R'_{\Box}\eta^*_{(S,\alpha,R_{\Box})}$.
Now we verify the condition 2 in the arrow part of Definition \ref{lrs}. Suppose $\eta^*_{(S,\alpha,R_{\Box})}(s)R'_{\Box}t$. Since $\eta^*_{(S,\alpha,R_{\Box})}$ is surjective, there is $t_1\in S$ such that $\eta^*_{(S,\alpha,R_{\Box})}(t_1)=t$ and $sR_{\Box}t_1$. Analogously we can verify for $\eta^{*}{-1}_{(S,\alpha,R_{\Box})}$. Therefore $\eta^*$ is a natural isomorphism.

\end{proof}
\end{thm}

\begin{figure}
\begin{center}
\includegraphics[width=100mm]{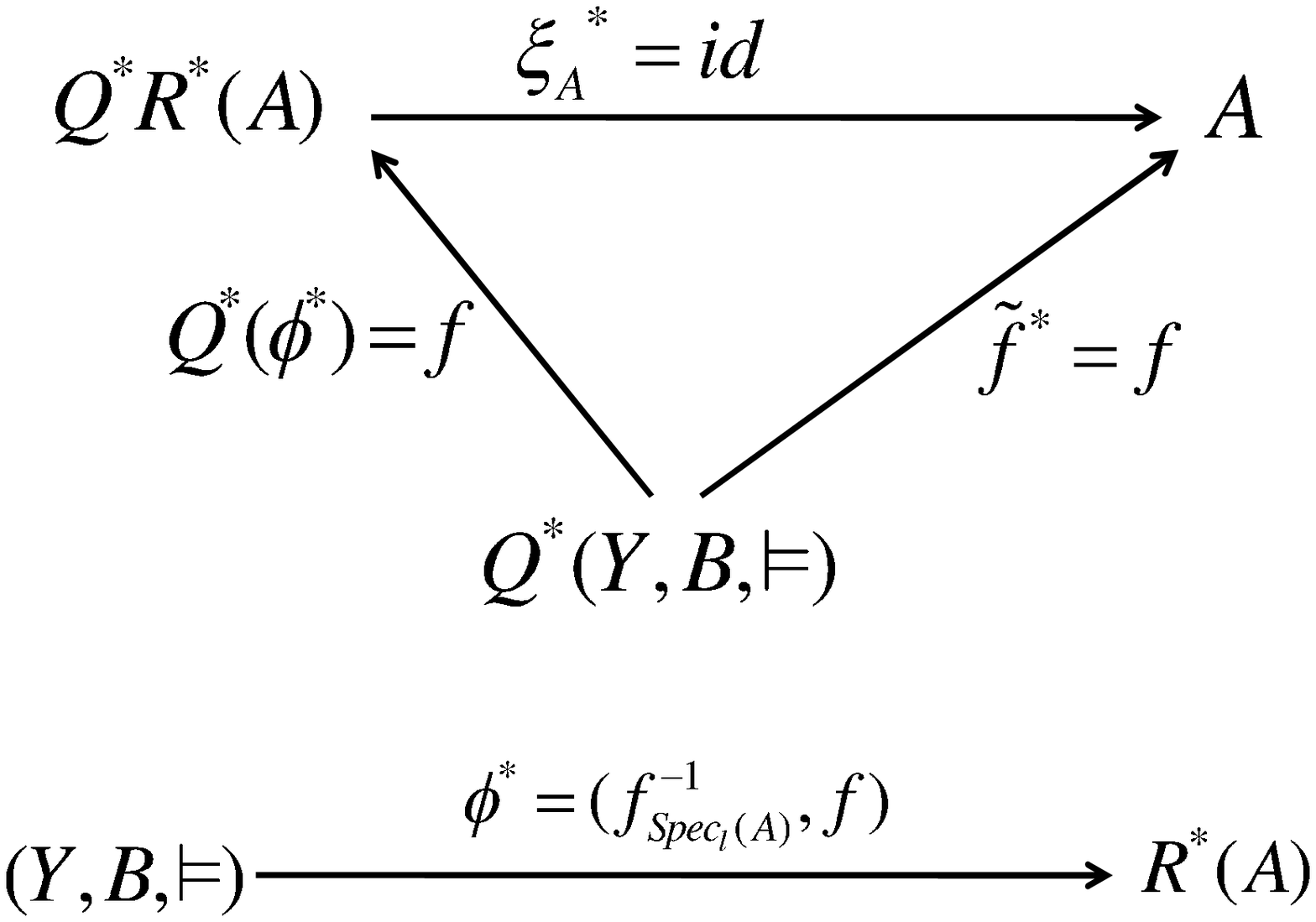}
\caption{ Diagram of Counit}
\label{fig:coadjoint222}
\end{center}
\end{figure}

\begin{thm}
\label{04}
$Q^*$ is the adjoint to the functor $R^*$.
\begin{proof}
We prove the theorem by unit of the adjunction. Diagram of unit is shown in Fig. \ref{fig:adjoint2}.

\begin{figure}
\begin{center}
\includegraphics[width=100mm]{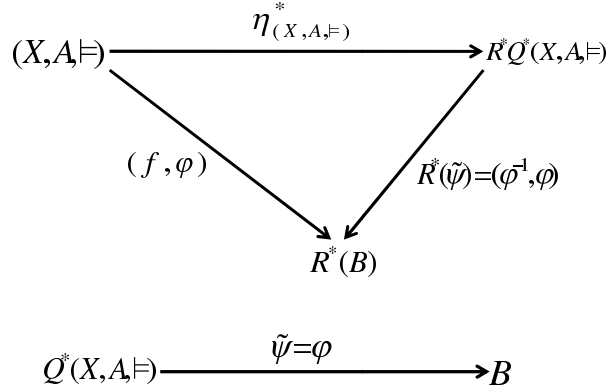}
\caption{Diagram of Unit}
\label{fig:adjoint2}
\end{center}
\end{figure}
Here $R^*(B)=(Spec_{\ell}(B), B, \models_{Spec_{\ell}\times B})$ where $\models_{Spec_{\ell}\times B}(v,b)=v(b)$. So $R^*Q^*(X,A,\models_{Spec_{\ell}\times B})=R^*(A)=(Spec_{\ell}(A), A, \models_{Spec_{\ell}(A)\times A})$.
Unit $\eta^*_{(X, A, \models)}: (X, A, \models)\rightarrow R^*Q^*(X, A, \models)$ is defined by $\eta^*_{(X, A, \models)}=(g,id_A)$.\\
Where
\begin{center}
\begin{enumerate}
\item $g: X\rightarrow Spec_{\ell}(A)$ is a set map. For each $x\in X$, we define $g(x)=g_x$. Where $g_x: A\rightarrow \ell$ such that $g_x(a)=\models(x,a)$.
\item $id_A: A\rightarrow A$ is an $\ell$-\textbf{ML}-algebra homomorphism.
\end{enumerate}
\end{center}
 It is already known that for each $x\in X$, $g_x$ is an $\ell$-\textbf{VL}-algebra homomorphism. From the Theorem \ref{a1} we observe that $(g,id_A)$ is a continuous map in $\ell$-\textbf{RSYM}. Now we shall show that the triangle given in the Fig \ref{fig:adjoint2} is commute i.e., 
 for a given arrow $(f,\varphi): (X,A,\models) \rightarrow R^*(B)$ there is an arrow $\tilde{\psi}$, we define $\tilde{\psi}=\varphi : Q^*(X, A , \models)\rightarrow B$ such that $(f,\varphi)= R^*(\tilde{\psi})\circ \eta^*_{(X,A,\models)}$. 
 Now $R^*(\tilde{\psi})= R^*(\varphi)=(\varphi^{-1},\varphi)$. So $f=\varphi^{-1}\circ g$ and $\varphi=id_A\circ \varphi$. Only part we have to show here is $f=\varphi^{-1}\circ g$.
  For each $s\in X$, $f(s)=\varphi^{-1}\circ g(s)=\varphi^{-1}\circ g_s= g_s\circ\phi$. Now for each $b\in B$, we have $ g_s\circ\varphi(b)=g_s(\varphi(b))=\models(s, \varphi(b))$.
  As $(f,\varphi)$ is continuous in $\ell$-\textbf{RSYM}, $\models(s,\varphi(b))=\models_{Spec_{\ell}\times B}(f(s),b)=f(s)(b)$. Therefore $f=\varphi^{-1}\circ g$. Hence the theorem is proved.

\end{proof}
\end{thm}
We can also prove $R^*$ is the coadjoint to the functor $Q^*$. Diagram of counit is shown in Fig. \ref{fig:coadjoint222}. 
For a given arrow $\tilde{f}^*$ in ($\ell$-MA)$^{op}$ there is an arrow in $\ell$-\textbf{RSYM}, we take $\phi^*=(f^{-1}_{Spec_{\ell}}(A), f)$ such that $Q^*(\phi^*)=f$. It can be readily shown that the triangle in the Fig. \ref{fig:coadjoint222} commute i.e., $\xi_A^*\circ Q^*(\phi^*)=\tilde{f}^*$.

\begin{thm}
\label{5}
The category ($\ell$-MA)$^{op}$ is equivallent to the category $\ell$-\textbf{RSYM}.
\begin{proof}
 We get two natural transformations $\xi^*$ and $\eta^*$ such that $\xi_A^*=id: Q^*R^*(A)\rightarrow A$ and $\eta^*_{(X,A,\models)}: (X,A,\models)\rightarrow R^*Q^*(X,A,\models)$. Here $\xi_A^*$ is conveniently a natural isomorphism. We have to show that $\eta^*_{(X,A,\models)}$ is a natural isomorphism between two lattice-relational systems. $\eta^*_{(X,A,\models)}=(g,id_A)$. Using the proof of the Theorem \ref{03}, we can say   $\eta^*_{(X,A,\models)}$ is a homeomorphism and hence $\eta^*_{(X,A,\models)}$ is a natural isomorphism. Therefore ($\ell$-MA)$^{op}$ is equivallent to the category $\ell$-\textbf{RSYM}. Consequently, $\ell$-MA is dually equivallent to the category $\ell$-\textbf{RSYM}.
 \end{proof}
 \end{thm}
 We now get the J\'{o}nsson-Tarski-style duality for $\ell$-\textbf{ML}-algebras. 
 \begin{thm}
 $\ell$-MA is dually equivallent to $\ell$-\textbf{RS}.
 \begin{proof}
 The result follows as a composition of the equivalences of Theorems \ref{3} and \ref{5}.
 \end{proof}
 \end{thm}
 \begin{rem}
The above duality which we already obtained in between $\ell$-MA and $\ell$-\textbf{RS} is the same as the the duality established in \cite{maruyama2011dualities}. We can easily apply the same procedure for proof of duality for $\ell$-\textbf{ML}-algebras with truth constants(c.f. \cite{maruyama2011dualities}).
\end{rem}

\section{Conclusion}
\label{sec:4}
 In this work, we study categorical relationships among $\ell$-\textbf{BSYM}, $\ell$-\textbf{BS} and $\ell$-VA-algebra. After this, we present the concept of lattice-relational systems and set up the equivalence with $\ell$-\textbf{RS}. Also we establish a duality between $\ell$-\textbf{ML}-algebra and system. Consequently, the duality between $\ell$-\textbf{RS} and $\ell$-\textbf{ML}-algebra has been worked out (In \cite{maruyama2011dualities}, the author has shown this duality directly). It has been found that during this progress of work, there were several parallel studies \cite{solovyov2012categorical,solovyov2010variable, rodabaugh1999categorical}, some of them introduce more generalized concept variety-based rahter than lattice (see \cite{solovyov2012categorical}) to present a categorical connection between system and space and some are taken fuzzy. In this regard, we may be noted that similar investigation can be continued by taking the algebraic model of any other modal logic instead of $\ell$-\textbf{ML}-algebra which is the algebraic model of Fitting's $\ell$-valued modal logic.\\
 In particular, we may provide below some future directions of the current set up.
 \begin{center}
 \begin{enumerate}
 \item Looking at the algebraic model of fuzzy(multi-valued) modal logic, we can move forward to similar research as revealed by our current work.
 \item By necessary developing of modal geometric logic, we may proceed for further generalization of the present paper.
 
 \end{enumerate}
 \end{center}

\end{document}